\documentclass[12pt,a4paper,twoside]{amsart}
\usepackage{amsmath}
\usepackage{amssymb}
\usepackage{amsfonts}
\usepackage{a4}
\usepackage{color}

\theoremstyle{plain}
\newtheorem{thm}{Theorem}
\newtheorem{cor}[thm]{Corollary}
\newtheorem{clm}[thm]{Claim}
\newtheorem{conj}[thm]{Conjecture}
\newtheorem{lem}[thm]{Lemma}
\newtheorem{lem*}[thm]{Lemma}

\theoremstyle{definition}
\newtheorem{dfn}{Definition}
\theoremstyle{remark}
\newtheorem{rem}{Remark}
\newtheorem{rem*}{Remark}

\newcommand{\mA}{\mathbb A}

\newcommand{\mC}{{\mathbb C}}

\newcommand{\mF}{\mathbb F}
\newcommand{\mG}{\mathbb G}

\newcommand{\mQ}{\mathbb Q}

\newcommand{\mZ}{{\mathbb Z}}

\newcommand{\ho}{\hookrightarrow}
\newcommand{\Gg}{\gamma}

\newcommand{\mcA}{\mathcal A}

\newcommand{\mcC}{\mathcal C}
\newcommand{\mcD}{\mathcal D}

\newcommand{\mcF}{\mathcal F}

\newcommand{\mcL}{\mathcal L}

\newcommand{\ti}{\tilde}

\begin{document}
\title{On the bias of cubic polynomials}
\begin{abstract}
Let $V$ be a vector space over a finite field $k=\mF _q$ of dimension $n$. For a polynomial $P:V\to k$ we define the {\it bias} of $P$ to be
$$b_1(P)=\frac {|\sum _{v\in V}\psi (P(V))|}{q^n}$$
where $\psi :k\to \mC ^\star$ is a non-trivial additive character.
 
A. Bhowmick and S. Lovett proved that for any $d\geq 1$ and $c>0$ there exists $r=r(d,c)$ such that 
any polynomial $P$ of degree $d$ with $b_1(P)\geq c$ can be written as a sum $P=\sum _{i=1}^rQ_iR_i$ 
where $Q_i,R_i:V\to k$ are non constant polynomials. We show the validity of a modified version of the converse statement for the case $d=3$.
\end{abstract}
\author{David Kazhdan, Tamar Ziegler}
\maketitle
\section{Introduction} 

 Let $k = \mF _q, q=p^m, k_n/k$ be the extension of degree $n$, $\bar k$ the algebraic closure of $k, tr_{k_n/k}:k_n\to k$ the trace. We fix a non-trivial additive character $\psi :k\to \mC ^\star$ and define $\psi _n:=\psi \circ tr_{k_n/k}:k_n\to \mC ^\star$. 

\begin{dfn}\label{1.1} Let $V$ be  a $k$-vector space and $PV\to k$
be a polynomial function  depending on a finite number of variables.   We define 

$(a)$ $a_n(P)=\sum _{v\in V(k_n)}\psi _n(P(v)),$

$(b)$ $\ti b_n(P)=\frac {|a_n(P)|}{ q^{ndim(V)}}, $

 $(c)$ $b_n(P)=-2\log _{q^n}(\ti b_n(P)) $
where $-log _{q^n}(0)=\infty$

and $(d)$ $b(P)=\liminf _nb_n(P)$.

We say that  $b_n(P)$ is the {\it bias} of $P$ over $k_n$ and $b(P)$ is the {\it bias} of $P$.
\end{dfn}
\begin{rem}  $a_n(V)$ is well defined since $a_n(P)=a_n(Q)$ if $P=f^\star (Q)$ where 
 $f:V\to W$ is a linear map. 
\end{rem}
\begin{dfn}\label{1.2} For a non-zero homogeneous  polynomial $P\in k[V^\vee ]$  of degree $d\geq 2$ we
define the {\it algebraic rank} $r(P)$ of $P$ as $2r$ where $r$ is the minimal number $r$ such that it is possible to write $P$ in the form
$$P = \sum ^ r_{i=1} l_iR_i$$ 
where $l_i, R_i\in \bar k[V^\vee ]$ are homogeneous polynomials of positive degrees. 
\end{dfn}

As shown in \cite{BL} for any $\geq 1,c>0$ there exists $M(d,c)$
such that $r(P)\leq M(d,c)$ 
 for any polynomial $P$ of degree $d$ with $b_1(P)\geq c$.

From now on we assume that $P$ is homogeneous of degree $d< p$.
\begin{thm}\label{1.3}   $b(P)<\infty$.
\end{thm}
\begin{cor}\label{1.4} For any  $l\geq 1, d<p$ there exist $n(l,d)\in \mZ _+, c(l,d)>0$ such that for any  polynomial $P:k^l\to k$ of degree $d$ we have $|a_n(P)|\geq c(l,d)$ for some  $n$, $1\leq n\leq n(l,d)$.
\end{cor}
The Corollary follows from the finiteness of $k$.

We want to express our gratitude to A. Braverman who helped with a proof of  Lemma \ref{4.3} and T. Schlank for a helpful discussions.

\section{ The quadratic  case } 

In this section we assume that $q$ is odd. Let $V$ be a finite-dimensional $k$-vector space, $P:V\to k$ a quadratic polynomial,

$P^\perp =\{v\in V|Q(v+v')=Q(v')$ for all $v'\in V$. Let $l:V\to k$ be a linear functional and 
$$a(P,l)=\sum _{v\in V(k)}\psi (P(v)+l(v))$$

\begin{lem}\label{2.1}

$(a)$ If $l_{P^\perp}\neq 0$ then $a(P,l)=0$.

$(b)$ If $l_{P^\perp}=0$ then $a(P,l)=0$ then $|r(P)-b(P)|\leq 1$. 

\end{lem}

\begin{prf} $(a)$. If  $l_{P^\perp}\neq 0$ then for any $v'\in V$ we have 
$$\sum _{v\in P^\perp (k)}\psi (P(v+v')+l(v+v'))=0.$$

$(b)$. If $l_{P^\perp}=0$ then there exist $v_0\in V,C\in k$ such that  
$$P(v)+l(v)=P(v+v_0)+C; \quad v\in V.$$
So we can assume that $l=0$. 

If $r=2t$ then there exists a system of coordinates $(x_1,...,x_N)$ on $V$ such that 
$$P(v)=\sum _1^tx_ix_{i+t}, v=(x_1,...,x_N)$$
Since $\sum _{x,y\in k}=\psi (xy)=q$ we have $b_1(P)=r$.

If $r=2t+1$ then there exists a system of coordinates $(x_1,...,x_N)$ on $V$ such that 
$$P(v)=\sum _1^tx_ix_{i+t}+\alpha x_r^2, \quad  \alpha \in k^\star,  v=(x_1,...,x_N)$$
In this case we have 
$b_1(P)=r-1/2$.

\end{prf}
\section{The cubic case.} 
\begin{thm}\label{3.1} For any $c>0$  there exist $t=t(c,q)>0$ such that 
 $b(P)\leq t(c,q)$ for  any cubic polynomial $P$  on a $k$-vector space $V$ with $r(P)\leq c$. Also there exists  $N(c,q)$ such 
 for  any cubic polynomial $P$  on a $k$-vector space $V$ with $r(P)\leq c$ we have $b_n(P)\leq t(c,q)$ for some $n\leq N(c,q)$.
\end{thm}

\begin{conj} One can choose $t(c)$ such that the statement of  Theorem \ref{3.1} is true for all finite fields of characteristic $>3$.
\end{conj}

In this section we derive  Theorem \ref{3.1} from Corollary \ref{1.4} and then prove Theorem \ref{1.3} in the next  section.

\begin{prf}

We have to show that for any $T\geq 1$ there exists $N(T),n(T)$ and $b(T)>0$ such that for any $k$-vector space $V$ and a cubic polynomial 
$P : V \to k$ such that $r(P)\leq T$ there exists $n\leq n(T)$ 
such that 
 $b_n(P)\geq b(T)$.

 We fix $P$ of the form   $P = \sum ^ r_{i=1} l_iR_i$ where $l_i$ are linear and $R_i$ quadratic polynomials. We can assume that the linear functions $l_i$ are linearly independent.

Let $W=\{ v\in V |l_i(v)=0,1\leq i\leq r\}$ and 
$R'_i$ be  the restrictions of $R_i$ on $W$.

Consider first the case when $r=1$.
We can assume that $V =k^N,l_1 (x_1,...,x_N)=x_1$. So 
$$P =x_1R(x)$$
 where $R$ is a quadratic form and $W=k^{N-1}$. We write $v\in V$ in the form  $v=(x,y),x\in k,y\in W(k)$.

For any  $x\in k^\star$ we define 
$$r_x=\sum_{y\in W(k)}\psi (R(x,y)).$$
By definition we have 
$$a(P)=\sum _{y\in W(k)}1+\sum _{x\in k^\star}r(a)=q^{(N-1)}+\sum _{x\in k^\star}r_a$$
It is sufficient to prove the following result.

\begin{lem}\label{3.2} Either  $|r_x|\leq q^{(N-3/2)}$ for all $x\in k^\star$ or 
there exist 
linear function $t_1,t_2,t_3$ on $V$ and a cubic form $Q$ on $k^3$ such that $P(v)=Q(t_1(v),t_2(v),t_3(v))$. 
\end{lem}


\begin{prf}
We write vectors of $v\in V$ in the form $v=(x,y),x\in k, y\in k^{N-1}$. Then
$$R(x,y)=\alpha x^2+xl(y)+Q(y)$$
 where $\alpha \in k, l$ is a linear function and $Q$ is a quadratic form on $W$. Let $W'\subset W$ be the kernel of $Q$. 

We consider two cases.

If $r(Q)\geq 3$ then it follows from Lemma \ref{2.1} that 
 $$|r_x|\leq q^{N-1} q^{-r(Q)/2}\leq q^{N-1-3/2}$$
 for all $x\in k^\star$ and therefore $b(P)=2$.

On the other hand if $r(Q)\leq 2$ then there is a subspace $W'\subset W$ of codimension $2$ such that $Q(y+y')=Q(y)$ for $y\in W,y'\in W'$. If $l_{W'}\neq 0$ then $r_x=0$ forall $x\in k^\star$.

One the other hand if $l_{W'}=0$ then $P=p^\star (\bar P)$ where $\bar P$ is a polynomial of degree $3$ on $W/W'$ and  Theorem \ref{3.1} follows from Corollary \ref{1.4}.
\end{prf}\\
 \ \\

To simplify notations we consider in details only the case $r=2$. The  general case is completely analogous. 

We can  assume that $V=k^N$ and $l_i(x_1,...,x_N)=x_i, 1\leq i\leq 2$. So 
$W\subset V$ is defined by the system of equation $\{ x_1(v)=x_2(v)=0\}$. We denote elements of $W$ by $y$ and write $v\in V$ in the form $v=(\bar x,y)$ where 
$\bar x=(x_1,x_2)$.
 Then we have 
$$P(\bar x,y)=\sum _{i=1}^2x_iR^i(y)$$
 where $R^i$ are polynomials of degree $2$. So there exist quadratic forms $c^i$ on $k^2$, linear maps $l^i\in Hom (k^2,W^\vee )$ and quadratic forms $Q^i$ on $W$ such that 
$$R^i(\bar x,y)=c^i(\bar x)+l^i(\bar x)(y)+Q^i(y), \bar x\in k^2,y\in W.$$
We define 
$$Q^{\bar x}=\sum _{i=1}^2 x_iQ^i$$
 for  $\bar x\in k^2-0$.

For any $\bar x\in k^2-0$ we define 
$$r_{\bar x}=\sum_{y\in W(k)}\psi (P(\bar x,y)).$$
By definition we have 
$$a(P)=q^{N-2}+\sum _{\bar a\in k^2-(0)}r_{\bar a}$$

If $rank (Q^{\bar x})\geq 12$ then (as follows from  Lemma \ref{2.1}) we have 
$$|r_{\bar x}|\leq q^{N-6}$$
Let $$U=\{ \bar x\in k^2-\{ 0\}|rank  (Q^{\bar x})\leq 12\}$$
and $R\subset k^2$ be the subspace generated by $U$. There are three possibilities for the dimension of $R$.

 $(0)$. If $R=\{ 0\}$ then 
$rank (Q^{\bar x})\geq 12$ for all $\bar x\in k^2-\{ 0\}$ we have $a_1(P)\sim q^{N-2}$.

$(2)$. If $R=k^2$ choose linearly independent $u,u'\in U$ and define 
$W'={Q^u}^\perp \cap {Q^{u'}}^\perp \subset W$. 
Then $W'\subset V$ is a subset of codimention $\leq 26$ and all quadratic forms $Q^{\bar x},\bar x\in k^2$ are invariant under shifts by $w'\in W'$. Let $W''=ker (l^1)\cap ker (l^2)\cap W'$.
Then there exists a polynomial $\ti P$ on $W/W''$ such that 
$P=\ti P\circ s$ where $s:W\to W/W''$ is the projection. It follows   from  Corollary \ref{1.4} that in this case
 the conclusion of Theorem   \ref{3.1} is true.

$(1)$. Now consider the case when $dim_k(R)=1$. We can then assume that 
$U=\{ (0,x)\},x\in k$.
Let  $V'=\{ (x_1,...,x_n)|x_1=0\}$ and $P'$ be the restriction of $P$ on $V'$. It is clear that there exists a cubic polynomial $\ti P'$ on $V'/W''$ such that  
$$P'=\ti P'\circ s', \quad s':V'\to V'/W''$$
Since $dim(V'/W'')\leq 26$ it  follows from  Corollary \ref{1.4} there exist $n\leq n(26,3)$ such that $|a_n(P')|\geq c$, $ c=c(26,3)$.

Let $u=-loq _q(c)$ and 
 $$U_1=\{ \bar x\in k^2-\{ 0\}|rank  (Q^{\bar x})\leq 16+2u\}$$
and $R'\subset k^2$ be the subspace generated by $U'$.

If $R_1=k^2$ then we finish the proof of Theorem \ref{1.4} as in the case $(2)$. On the other hand if $dim(R_1)=R$ we see that 
$|a_n(P)-a_n(P')|\leq q^{(N-u-1)n}$ and therefore as before 
it follows   from the definition of $u$ that in this case
 the conclusion of  Theorem  \ref{3.1} is also true.

We only indicate a way  to start  the  proof of the result for the general $r$ and leave details to the reader.

We can  assume that $V=k^N, P=\sum _{i=1}^rx_i\ti Q^i$ where $\ti Q^i$ are quadratic forms. Let 
$W\subset V$ is defined by the system of equation $\{ x_1(v)=...=x_r(v)=0\}$. We denote elements of $W$ by $y$ and write $v\in V$ in the form $v=(\bar x,y)$ where 
$\bar x=(x_1,...,x_r)$.
 Then we have 
$$P(\bar x,y)=\sum _{i=1}^rx_iR^i(y)$$
where $R^i$ are polynomials of degree $2$. So there exist quadratic forms $c^i$ on $k^r$, linear maps $l^i\in Hom (k^r,W^\vee )$ and quadratic forms $Q^i$ on $W$ such that 
$$R^i(\bar x,y)=c^i(\bar x)+l^i(\bar x)(y)+Q^i(y), \bar x\in k^r,y\in W$$
We define 
$$Q^{\bar x}=\sum _{i=1}^r x_iQ^i$$
 for  $\bar x\in k^r-0$.

For any $\bar x\in k^r-0$ we define 
$$r_{\bar x}=\sum_{y\in W(k)}\psi (P(\bar x,y)).$$
By definition we have 
$$a(P)=q^{N-r}+\sum _{\bar a\in k^r-(0)}r_{\bar a}$$
Let $u_0=-loq _qc(2r+2,3)$ and define
$$U_0:=\{ \bar x\in k^r-\{ 0\}|rank  (Q^{\bar x})\leq 2r+2+u_0\}$$
and denote by $R_0\subset k^r$  the subspace generated by  $U_0$.
If $U_0=\{ 0\}$ then as before we find that 
$a_1(P)\sim q^{N-r}$ and we see that  in this case
 the conclusion of  Theorem  \ref{3.1} is true.

If $R_0\neq \{ 0\}$ then we define $U_1$ as earlier in our proof of the case when $r=2$.

\end{prf}
\section{A proof of Theorem \ref{1.3}} 
We start notations from section $1.1$ in \cite{D}. For any algebraic $k$-variety 
$X$ we denote by  $\mcD (X)$ be the bounded derived category of complexes $\mF$ of constructible $l$-adic sheaves on $X$. For any $F,G$ in  $\mcD (X)$ we denote by $F\otimes G\in Ob(\mcD (X))$  the derived tensor product. 
A morphism $p:X\to Y$ of algebraic varieties 
defines (derived) functors $p_\star ,p_!
:\mcD (X)\to \mcD (Y)$
and  $f^\star  :\mcD (Y)\to \mcD (X)$.

In the case when $Y=\star =Spec(k)$ we can identify cohomology of compexes  $p_\star (F)$ and $p_!(F)$ 
 of $\mQ _l$-vector spaces with the cohomology groups 
$H^\star (X,F)$ and $H_c^\star (X,F)$. On the other hand for any point $x\in X(k)$ the stalk $F_x$ of $F$ at $x$ is equal to $i_x^\star (F)$ which is an object of $\mcD (Spec (k))$ where 
 $i_x:Spec (k)\to X$ is the imbedding $x\ho X$. We denote by 
$E_x(F)$ the Euler characteristic of the complex $i_x^\star (F)$. 

Let $\psi$ be a non-trivial additive character on $k$ with values in $\bar \mQ _l$ where $l$ is prime number $(q,l)=1$.
Deligne defined the {\it Fourier transform} functor $\mcF \in Aut (\mcD)(\mA )$ by 
$$\mcF (F)={p_1}_! (p^\star (\mcL _{AS}\otimes p_2^\star (F))$$
where $F$ is an object of $\mcC, p_1,p_2,p:\mA ^2\to \mA$ 
maps given by 
$$p(x,y)=xy, p_1(x,y)=x, p_2(x,y)=y,$$
(see the section $3$ of \cite{La}) where $\mA$  is the affine line over $k$.

The following result follows from results of \cite{La}.

\begin{clm} If  $F\in Ob(\mcD (\mA ))$ is $\mG _m$-equivariant then 
the Euler characteristic of the stalk $\mcF (F)_1$ at $1\in \mA (k)=k$ is equal to the difference  $E_0(\mcF (F))-E_1(\mcF (F))$.
\end{clm}

One can identify the set  of $\mA (k)$ with $k$. 
We define the {\it Frobenious map} $Fr_\mA :\mA \to \mA$ by $x\to x^q$. Then $k$ is the set of fixed points of $Fr_\mA$.  For any $n\geq 1$ we define $\psi _n:k_n\to \bar \mQ _l$ by $\psi _n=\psi \circ tr (k_n/k)$ where $k_n=\mF _{q^n}$.

Let $ p:\mA \to \mA$ be the map given by $x\to x^q-x$. Then $p:\mA \to \mA$ is a Galois covering with the group equal to the additive group of $k$ which acts by 
 $$\Gg :x\to x+\Gg  ,\quad \Gg \in k, x\in \mA.$$

 It defines 
a {\it local system}  $\mcL ^{AS}$  on $\mA$ such that fibers $\mcL ^{AS}_x$ of $\mcL$ at $x\in \mA$  are one-dimensional $\bar \mQ _l$-vector spaces. Any morphism 
$$\alpha :Fr_\mA ^\star (\mcL ^{AS}) \to \mcL ^{AS}$$
defines a map 
$$\alpha _x:\mcL _{x^q}\to \mcL ^{AS}_x$$
for any $x\in \mA$. In particular for any $x\in k, \alpha _x$ 
is an endomorphism of the line $\mcL _x$. So we can consider it 
as an element of $\bar \mQ _l$.

 \begin{clm}\label{4.1}  There exists unique  isomorphism 
$$\alpha :Fr_\mA ^\star (\mcL) \to \mcL$$
such that $\alpha _x=\psi (x), x\in k$.
\end{clm}
Now let $P:\mA ^N\to \mA$ be a polynomial defined over $k, \mcL _P=P^\star (\mcL ^{AS})$. Then $\mcL _P$  is a local system on $\mA ^N$ and $\alpha$ defines an isomorphism 
$$\beta :Fr_{\mA ^N}^\star (\mcL _P) \to \mcL _P$$
such that 
$$\beta _v=\psi (P(v)), v\in k^N\subset \mA ^N$$
Let $F=p_!(\mcL _P), G=\mcF (F)$. 
$\beta$ induces  morphisms $\Gg:Fr_\mA ^\star (F)\to F$ and 
$\ti \Gg :Fr_\mA ^\star (G)\to G$.
We denote by $\Gg _a, \ti \Gg _a$ the induced endomorphisms of complexes $F_a,G_a, a\in k$.
For any $i$ we denote by $\beta _i$ the induced endomorphism of 
$H^i_c(\mcA ^N,\mcL _P)$.

The following equality follows from the Lefschetz fixed-point theorem 
\begin{clm}\label{4.2}  For any $n\geq 1$ we have 
 $$a_n(P)=\sum _i(-1)^itr(\beta _i^n)=Tr (\ti \Gg )_1^n$$
\end{clm}

We define 
$$E(P):=a_0(P)=\sum _i(-1)^idim (H^i_c(\mA ^N,\mcL _P)) $$
As follows from this Claim we have $E(P)=E(G_1)$.

\begin{lem}\label{4.3}  If  $P$ is homogeneous polynomial of degree $d$ prime to $p$ then $E(P)\neq 0$.
\end{lem}

\begin{prf} 
Since $P$ is  homogeneous the object $F$ of $\mcD (\mA)$ is $\mG _m$-equivariant. So the Claims \ref {4.1} ,\ref{4.2} imply that

$$E(P)=E(F_0)-E(F_1)=E(X_0)-E(X_1)$$ 
where $X_a:=P^{-1}(a)\subset \mA ^N$ and 
$$E(X):=\sum _i(-1)^idim (H^i_c(X),\bar \mQ _l))$$
for any algebraic $k$-variety $X$.

Since the group $\mG _m$ acts freely on $X_0-\{ 0\}$ we see that 
$E(X_0)=E(\{ 0\}) =1$. On the other hand the group $\mu _d$ of $d$-roots of unity acts freely on $\mcA ^N-\{ 0\}$ by 
$$\xi (x_1,...,x_N)=(\xi x_1,...,\xi x_N), \xi ^d=1$$preserving $X_1\subset \mA ^N$. Let $Y=E(X_1)/\mu _d$.
 
Then  $E(X_1)=dE(Y)\equiv 0$ mod $d$. So $E(P)\neq 0$.
\end{prf}

Theorem \ref{1.3} follows now from Claim   \ref{4.2}.

\section{A  very special case, } 
 In the case when $P=Q^2$ where $Q$ is a non-degenerate quadratic form on $V=k^N$ we have $|a_n(P)|\sim q^{nN/2}$. We see that for quartic polynomials $P$  we can not bound $b(P)$ below in terms of $r(P)$. On the other hand such a bound exists in {\it nice} cases.
\begin{lem}{5.1} 
Assume that  $P=QR:\mA ^N\to \mA, deg(Q),deg(R)\neq 0$,  that all  irreducible components of 
$$X=\{ v\in \mA ^N|P(v)=0\}$$
 are reducible and defined over $k$ and that $(d,q-1)=1$. There exists a constant $\Gg (d)$  such that
$$|a_1(P)|\geq q-\Gg (d)\sqrt q$$

\end{lem}
\begin{prf}  Since $(d,q-1)=1$ all the fibers $P^{-1}(a)(k),a\in k^\star$ have the same number of elements $A$. It is clear that 
$A\sim q^{N-1}$. Let $B=|P^{-1}(0)(k)|=|X(k)|$. Then $a_1(P)=B-A$.
On the other hand $X=Y\cup Z$ where
$$Y=\{ v\in \mA ^N|Q(v)=0\}, Z=\{ v\in \mA ^N|R(v)=0\}$$
 So $B=|Y(k)|+|Z(k)|-|Y\cup Z(k)|$ and the Lemma follows from the Weil's estimates (see \cite{D}).
\end{prf}


\begin{thebibliography}{99}
\bibitem[D]{D} 
Deligne, Pierre
{\em La conjecture de Weil. II},
Inst. Hautes ƒtudes Sci. Publ. Math. No. 52 (1980), 137Ð252. 

	\bibitem[BL]{BL} Abhishek Bhowmick, Shachar Lovett {\em Bias vs structure of polynomials in large fields, and applications in effective algebraic geometry and coding theory}. 


\bibitem[La]{La} Laumon, G. {\ Transformation de Fourier, constantes d'Žquations fonctionnelles et conjecture de Weil.} 
Inst. Hautes ƒtudes Sci. Publ. Math. No. 65 (1987), 131Ð210. 
\bibitem[L]{L} Laumon, G. {\em Exponential sums and l-adic cohomology: a survey}, Israel J. Math. 120 (2000), part A, 225Ð257. 


\end{thebibliography}
\end{document}